\documentclass[titlepage,11pt]{article}
\usepackage{latexsym}
\usepackage{amsfonts}
\usepackage{amsmath}
\usepackage{tikz}
\usepackage{float}
\usepackage{url}
\usetikzlibrary{arrows}
\usetikzlibrary{decorations.markings}
\newcommand\blackslug{\hbox{\hskip 1pt \vrule width 4pt height 8pt depth 1.5pt
        \hskip 1pt}}
\newcommand\bbox{\hfill \quad \blackslug \bigbreak}

\title{Proof of the Caccetta-H\"{a}ggkvist conjecture for digraphs with small independence number}
\author{Patrick Hompe}

\date{\today}

\newtheorem{thm}{Theorem}[section]
\newtheorem{conjecture}{Conjecture}[section]
\newtheorem{lemma}{Lemma}[section]

\newcommand{\Proof}{\noindent{\bf Proof.}\ \ }

\begin{document}
\maketitle
\begin{abstract}
For a digraph $G$ and $v \in V(G)$, let $\delta^+(v)$ be the number of out-neighbors of $v$ in $G$. The Caccetta-H\"{a}ggkvist conjecture states that for all $k \ge 1$, if $G$ is a digraph with $n = |V(G)|$ such that $\delta^+(v) \ge n/k$ for all $v \in V(G)$, then G contains a directed cycle of length at most $k$. In [2], N. Lichiardopol proved that this conjecture is true for digraphs with independence number equal to two. In this paper, we generalize that result, proving that the conjecture is true for digraphs with independence number at most $(k+1)/2$.
\end{abstract}

\section{Introduction and definitions}
For the rest of the paper, we use the words cycle and path to refer to a directed cycle and directed path, respectively, and every graph considered is a digraph. Furthermore, every digraph $G$ is simple, meaning it has no loops or parallel edges. Let the girth $g(G)$ of a digraph $G$ be the length of its shortest cycle, and for a vertex $v \in V(G)$, let $\delta^+(v)$ denote the number of out-neighbors of $v$ in $G$. Let $\Delta^+(G) = \displaystyle \min_{v \in V(G)} \delta^+(v)$ be the minimum out-degree of a vertex in $G$. For vertices $u,v \in V(G)$, let the distance $d(u,v)$ from $u$ to $v$ be the length of the shortest path from $u$ to $v$ (define this to be zero if $u = v$). For $v \in V(G)$ and $i \ge 1$, let $N_i^+(v)$ be the set of vertices $u$ with $d(v,u) = i$, and let $N_i^-(v)$ be the set of vertices $u$ with $d(u,v) = i$. For a digraph $G$, call a set of vertices $H \subset V(G)$ independent if there are no edges between any two vertices of $H$. Let the independence number $\alpha(G)$ of a digraph $G$ be the size of the largest independent set $H \subset V(G)$. For disjoint sets $S_1, S_2 \subset V(G)$, say that $S_1$ is stable with $S_2$ if there are no edges between a vertex in $S_1$ and a vertex in $S_2$.

We begin with the following simple observation.
\begin{lemma}\label{trivial}
Suppose that $G$ is a digraph containing a cycle; then $g(G) \le 2\alpha(G)+1$.
\end{lemma}
\Proof
Let $C$ be a cycle of $G$ with minimum length, and suppose $C$ has at least $2\alpha(G)+2$ vertices. Then there exists a subset $S \subset V(C)$ of size $\alpha(G)+1$ such that no pair of vertices of $S$ are adjacent in $C$. Then there is an edge in $G$ between some pair of vertices in $S$, which gives a shorter cycle in $G$, a contradiction. This proves Lemma \ref{trivial}.\bbox{}

The next lemma immediately follows from Lemma \ref{trivial}, and is used repeatedly throughout the paper.

\begin{lemma}\label{acyclic}
Suppose $G$ is a digraph with $g(G) \ge 2\alpha(G)$, and that $H \subset G$ is a subgraph of $G$ with $\alpha(H) \le \alpha(G) - 1$. Then $H$ is acyclic.
\end{lemma}
\Proof
If $H$ contains a cycle, then Lemma \ref{trivial} shows that $H$ contains a cycle of length at most $2\alpha(G)-1$, which is a contradiction. This proves Lemma \ref{acyclic}.\bbox{}

In this paper, we deal with the following formulation of the Caccetta-H\"{a}ggkvist conjecture, which was introduced in [\ref{ch}]:
\begin{conjecture}[Caccetta-Haggkvist]\label{CH}
For $d \ge 1$, $k \ge 1$, if $G$ is a digraph with $n = |V(G)| \le kd$ and $\Delta^+(G) \ge d$, then $g(G) \le k$.
\end{conjecture}

For $k=1$ and $k=2$ it follows that the digraph is not simple, a contradiction. So, to prove Conjecture \ref{CH}, we can assume $k \ge 3$.

Now, Lemma \ref{trivial} gives that Conjecture \ref{CH} is true for $\alpha(G) \le (k-1)/2$. In this paper, we prove that Conjecture \ref{CH} is true for $\alpha(G) \le (k+1)/2$.

\section{Main Results}
We need the following two lemmas.
\begin{lemma}\label{lemma1}
Suppose that $G$ is an acyclic digraph; then for all $v \in V(G)$, there exists a path of length at most $2\alpha(G)-1$ to a vertex of out-degree zero in $G$.
\end{lemma}
\Proof
Since $G$ is acyclic, there exists a path from $v$ to a vertex of out-degree zero in $G$. Let $P = (v, v_2, \cdots, v_k)$ be a shortest such path. Then $P$ is induced, so if $k \ge 2\alpha(G) + 1$ then $\{v, v_3, \cdots, v_k\} \subset V(G)$ is an independent set of size at least $\alpha(G) + 1$, which is a contradiction. Thus $P$ has length at most $2\alpha(G) - 1$, as desired. This proves Lemma \ref{lemma1}.\bbox{}
\begin{lemma}\label{lemma2}
Let $G$ be an simple digraph with minimum out-degree $d \ge 1$, $\alpha(G) \ge 3$, and $g(G) \ge 2\alpha(G)$. Set $p = 2\alpha(G)-3$, and suppose $v \in V(G)$ is a vertex with $\delta^+(v) = d$. For odd $1 \le i \le p$, let $S_i$ be the subgraph of $G$ induced by the vertex set $V(G) \setminus (N_1^+(v) \cup \{v\} \cup \left(\bigcup_{j=1}^i N_j^-(v)\right)$. Then, for each odd $1 \le i \le p$, there exists a unique vertex $v_i \in S_i$ such that $N_1^+(v_i) \subset N_i^-(v)$. Furthermore, $|V(G) \setminus S_p| \ge (2\alpha(G)-2)d+1$.
\end{lemma}
\Proof
Every $w \in N_1^+(v)$ has $N_1^+(w) \subset N_1^+(v) \cup S_p$, since otherwise we obtain a cycle of length at most $2\alpha(G)-1$, a contradiction. Since $|N_1^+(v)|=d$, it follows that $V(S_i) \ne \emptyset$ for odd $1 \le i \le p$.

Now, for odd $1 \le i \le p$, we iteratively choose $v_i \in S_i$ such that $N_1^+(v_i) \subset N_i^-(v)$. Let $\{v_1,v_3,\cdots,v_{i-2}\}$ be vertices such that $N_1^+(v_j) \subset N_j^-(v)$ for all odd $1 \le j \le i-2$ (if $i=1$, this set of vertices is empty). The set $T = \{v, v_1,v_3,\cdots,v_{i-2}\}$ (if $i=1$, then $T = \{v\}$) is stable with $S_i$, so $\alpha(S_i) \le \alpha(G)-(i+1)/2$, and thus $S_i$ is acyclic by Lemma $\ref{trivial}$. Thus there exists $v_i \in S_i$ with out-degree zero in $S_i$.

Now, we claim that $N_1^+(v_i) \subset N_i^-(v)$. If not, then $v_i$ has an edge to a vertex $w_1 \in N_1^+(v)$, which has an edge to a vertex $w_2 \in S_i$. Let $H$ be the subgraph of $S_i$ induced by the set of vertices with no edge to $v_i$. We may assume $w_2 \in H$. We have that $\{v, v_1, \cdots, v_i\}$ is stable with $H$, so $\alpha(H) \le \alpha(G) - (i+3)/2$. Then by Lemma $2.1$ there exists a path $(w_2 \cdots w_j)$ of length at most $2\alpha(G)-i-4$ from $w_2$ to a vertex $w_j \in H$ with out-degree zero in $H$. If $w_j$ has out-degree in $S_i$ equal to zero, then since $|N_1^+(v)| = d$, it follows that $w_j$ has an out-neighbor in $N_i^-(v)$ and we obtain a cycle of length at most $2\alpha(G)-1$, a contradiction. If instead $w_j$ has an out-neighbor to $w_{j+1} \in S_i \setminus H$, then we again obtain a cycle of length at most $2\alpha(G)-1$, a contradiction. It follows that $v_i$ has $N_1^+(v_i) \subset N_i^-(v)$ for odd $1 \le i \le p$, as claimed.

Now, for odd $1 \le i \le p$, let $V_i$ be the set of vertices $u \in S_i$ such that $u$ has out-degree zero in $S_i$. Let $H = \{v\} \cup V_1 \cup V_3 \cup \cdots \cup V_p$. For $v_i \in V_i$, since $N_1^+(v_i) \subset N_i^-(v)$, it follows that $H$ is an independent set, so $|V(H)| \le \alpha(G)$. We also know that the $V_i$ are nonempty, so $|V(H)| \ge \alpha(G)$. Thus $|V_i| = 1$ for all odd $1 \le i \le p$. This proves the first part of the lemma, namely that for each odd $1 \le i \le p$ there exists a unique vertex $v_i \in S_i$ with $N_1^+(v_i) \subset N_i^-(v)$. For the remainder of the proof, let $\{v_1,v_3,\cdots,v_p\}$ be those unique vertices.

For odd $3 \le i \le p$, define $X_i = N_1^+(v_i) \subset N_i^-(v)$, and let $T_i = N_i^-(v) \cup N_{i-1}^-(v) \setminus X_i$. $\{v\}$ is stable with $X_i$, so by Lemma \ref{acyclic}, $X_i$ is acyclic and contains a vertex $u_i \in V(X_i)$ with out-degree zero in $X_i$. We claim that $N_1^+(u_i) \subset T_i$, and consequently $|T_i| \ge d$. If not, then there exists a path of length at most 2 from $u_i$ to a vertex $w_2 \in S_i$. Since $\{v, v_1, \cdots v_{i-2}\}$ is stable with $S_i$, $\alpha(S_i) \le \alpha(G) - (i + 1)/2$. Lemma \ref{lemma1} gives a path from $w_2$ to $v_i$ of length at most $2\alpha(G) - 5$. These two paths together form a cycle in $G$ of length at most $2\alpha(G) - 2$, which is a contradiction.

Thus, for odd $3 \le i \le p$, we have $|X_i| + |T_i| \ge 2d$. Also, $N_1^+(v_i) \subset N_1^-(v)$ gives $|N_1^-(v)| \ge d$, and by the definition of $v$ we have $|N_1^+(v)| = d$. Together with the vertex $v$, these give:
\[|V(G) \setminus S_p| \ge (p-1)d + 2d + 1 = (2\alpha(G)-2)d+1\]
as desired. This proves Lemma \ref{lemma2}.~\bbox{}

Lemma \ref{lemma2} is used to prove the following two theorems.

\begin{thm}\label{2alpha}
Suppose that $G$ is a digraph with minimum out-degree $d \ge 1$ and $n = |V(G)| \le 2\alpha(G)d$; then $g(G) \le 2\alpha(G)$.
\end{thm}
\Proof
As mentioned above, it suffices to consider simple digraphs $G$ with $\alpha(G) \ge 2$. The case $\alpha(G) = 2$ is proved in [\ref{lich}], so we may further assume that $\alpha(G) \ge 3$. Now, for the sake of contradiction, suppose that $g(G) \ge 2\alpha(G) + 1$. Then Lemma \ref{lemma2} implies that $|V(G) \setminus S_p| \ge (2\alpha(G)-2)d+1$, which together with $|V(G)| \le 2\alpha(G)d$ gives $|S_p| \le 2d-1$. $H=\{v, v_1, v_3, \cdots, v_{p-2}\}$ is stable with $S_p$, so $\alpha(S_p) = 1$ and $S_p$ is a transitive tournament. Let $(w_1\cdots w_r)$ be the unique Hamiltonian path of the transitive tournament $S_p$.

Now, $J = \{v_1,v_3,\cdots,v_p\}$ is stable with $N_1^+(v)$, so $N_1^+(v)$ is a transitive tournament. Let its unique Hamiltonian path be $(u_1\cdots u_d)$. $N_1^+(u_d) \subset S_p$, so there is an out-neighbor $w_k$ of $u_d$ with $k \ge d$. It follows that $w_k$ has an edge to a vertex not in $S_p$. An edge from $w_k$ to $v$ or to $w' \in N_i^-(v)$ for some $1 \le i \le p$ yields a cycle of length at most $2\alpha(G)$, a contradiction. If, instead, $w_k$ has an edge to $u' \in N_1^+(v)$, then $u'$ has an edge to $u_d$, and we obtain a cycle of length at most three, a contradiction. This proves Theorem \ref{2alpha}.~\bbox{}

\begin{thm}\label{2alpha-1}
Suppose $G$ is a digraph with minimum out-degree $d \ge 1$ and $n = |V(G)| \le (2\alpha(G)-1)d$; then $g(G) \le 2\alpha(G)-1$.
\end{thm}
\Proof
As mentioned above, it suffices to consider simple digraphs $G$ with $\alpha(G) \ge 2$. The case $\alpha(G) = 2$ is proved in [\ref{lich}], so we further assume that $\alpha(G) \ge 3$. For the sake of contradiction, suppose $g(G) \ge 2\alpha(G)$. Lemma \ref{lemma2}
gives a set of vertices $\{v_i\}$ indexed by odd $1 \le i \le p$ such that $N_1^+(v_i) \subset N_i^-(v)$. $v_1$ is stable with $N_1^+(v)$ (otherwise we obtain a cycle of length four), so Lemma \ref{acyclic} gives that $N_1^+(v)$ is acyclic. So, there exists $u \in N_1^+(v)$ with out-degree zero in $N_1^+(v)$. If $u$ has an edge to any vertex not in $S_p$, we obtain a cycle of length at most $2\alpha(G)-1$, a contradiction. Thus, we must have $N_1^+(u) \subset S_p$ and it follows that $|S_p| \ge d$.

But Lemma \ref{lemma2} also gives that $|V(G) \setminus S_p| \ge (2\alpha(G)-2)d+1$, which together with $|S_p| \ge d$ implies that $|V(G)| \ge (2\alpha(G)-1)d+1$, contradicting the assumption that $|V(G)| \le (2\alpha(G)-1)d$. This proves Theorem \ref{2alpha-1}.~\bbox{}

\begin{thm}\label{mainresult}
Conjecture \ref{CH} is true for digraphs $G$ with $\alpha(G) \le (k+1)/2$.
\end{thm}
\Proof
Theorem \ref{2alpha} and Theorem \ref{2alpha-1} together with Lemma \ref{trivial} give the desired result. This proves Theorem \ref{mainresult}. ~\bbox{}

\end{document}